\documentclass[11pt]{article}
\usepackage{arxiv}

\usepackage{times}
\usepackage{graphicx}        
\usepackage{multicol}        
\usepackage[bottom]{footmisc}

\usepackage{amssymb,amsmath,mathrsfs,amsthm}
\usepackage{mathtools}
\let\labelindent\relax
\usepackage[inline,shortlabels]{enumitem} 
\usepackage{framed}
\usepackage{algorithm,algorithmic,psfrag}
\usepackage{soul}
\usepackage{color}
\usepackage{kbordermatrix}
\usepackage[utf8]{inputenc} 
\usepackage[T1]{fontenc}    
\usepackage{hyperref}       
\usepackage{url}            
\usepackage{booktabs}       
\usepackage{amsfonts}       
\usepackage{nicefrac}       
\usepackage{microtype}      
\usepackage{doi}
\usepackage{xifthen}

\usepackage{tikz}
\usetikzlibrary{positioning,calc}
\usetikzlibrary{fit}

\newtheorem{theorem}{Theorem}
\newtheorem{corollary}{Corollary}

\definecolor{paleGreen}{rgb}{.3, .7, .3}
\definecolor{coolBlue}{rgb}{.3, .5, 1}
\definecolor{rosePink}{rgb}{.9, .5, .4}

\newcommand{\DpCoeff}[1]{\alpha_{#1}}
\newcommand{\DqCoeffpr}[1]{\beta_{#1}}
\newcommand{\DqCoeffpl}[1]{\kappa_{#1}}
\newcommand{\DqCoeffql}[1]{\gamma_{#1}}

\newcommand{\R}{\mathbb{R}}

\newcommand{\zero}{\mathbf{0}}
\newcommand{\one}{\mathbf{1}}
\newcommand{\eye}{I}
\DeclareMathOperator*{\diag}{diag}
\newtheorem{defin}{Definition}

\title{Control-Oriented Modeling of Pipe Flow through Intersecting Pipe Geometries}

\author{ \href{https://orcid.org/0000-0001-8308-4331}{\includegraphics[scale=0.06]{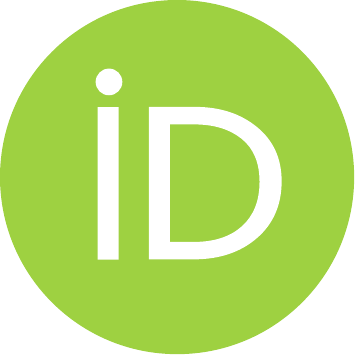}\hspace{1mm}Sven Br{\"u}ggemann} \\
	Mechanical \&\ Aerospace Engineering Department\\
	University of California, San Diego\\
	CA 92093-0411, USA\\
	\texttt{sbruegge@eng.ucsd.edu} \\
	\And
	\href{https://orcid.org/0000-0002-5067-1940}{\includegraphics[scale=0.06]{orcid.pdf}\hspace{1mm}Robert R.~Bitmead} \\
Mechanical \&\ Aerospace Engineering Department\\
	University of California, San Diego\\
	CA 92093-0411, USA\\
	\texttt{rbitmead@eng.ucsd.edu} \\
}

\begin{document}
\maketitle
\begin{abstract}
	We present control-oriented models for transient dynamics of isothermal one-dimensional gas flow through multiple pipes in series and intersecting pipe geometries. These composite models subsume algebraic constraints that would otherwise appear due to boundary conditions, so that our linear state-space models are well-suited for model-based control design for gas flow in pipe networks with non-trivial geometries.
\end{abstract}
\section{Introduction}
For transient modeling of gas flow through pipe networks, fluid dynamics and, particularly, computational fluid dynamics, are well-established subjects focused on high-fidelity modeling given design and boundary conditions; typically, they involve nonlinear partial differential equations (PDEs) and transport phenomena which are not amenable to finite-dimensional control design but instead are targeted and tested for simulation.

Other more pragmatic modeling for gas pipeline distribution systems \cite{benner2019,kralik1988,Behbahaninejad2008AMS} is usually based on discretization and  yields a system of ordinary differential algebraic equations (DAEs), which again is not well-suited to control design. Although, it can be used directly for controller synthesis in some circumstances \cite{mpdDae} and, as noted in \cite{benner2019}, if the DAE is of index~1. Theorem~4.1 in Benner \text{et al.} \cite{benner2019} establishes that the DAEs describing the gas flow through interconnected pipes are indeed of index~1. 

This fact is used in \cite{sven_bob_rob_gas} to rewrite the system of DAEs as a (state-space) system of linear ordinary differential equations (ODEs) \textit{subsuming} the algebraic constraints. Thus, to synthesize model-based controllers the rich literature on Linear Systems Theory can be exploited such as the Mason's Gain Formula for modeling interconnections. An equivalent approach modified for state-space realizations is presented in \cite{sven_bob_rob_gas}.

Aiming for control-oriented models of non-trivial pipe geometries akin to \cite{sven_bob_rob_gas} the purpose of this technical report is twofold. Firstly, we linearize the nonlinear ODEs characterizing an isothermal one-dimensional gas flow through a single pipe from \cite{benner2019} and recall the related model for a flow through multiple pipes in series from \cite{sven_bob_rob_gas}. We provide a brief analysis in the frequency domain. Secondly, again based on the single-pipe dynamics, we extend the models for joining and branching pipe flows from \cite{sven_bob_rob_gas} to non-trivial pipe geometries of arbitrary number of pipes.

All the aforementioned models appeal to model-based MIMO control as they are in state-space form and they can be parametrized by physics. 
\section{Preliminaries: single-pipe model}
Under the assumptions of a constant temperature and a flow velocity much lower than the speed of sound Benner \emph{et al.} \cite{benner2019} derive a nonlinear model for an isothermal one-dimensional pipe flow. Towards a discretization and linearization of the related nonlinear dynamics, denote $p(x,t)$ as the pressure and $q(x,t)$ as the mass flow. Let the boundary conditions
\begin{align*}
p_\ell\doteq p(0,t), \qquad q_r\doteq q(L,t)
\end{align*}
be given, whereas 
\begin{align*}
p_r\doteq p(L,t), \qquad q_\ell\doteq q(0,t)
\end{align*}
are to be determined through the model.
Spatial discretization and linearization of \cite[Eq. 3.2]{benner2019} yields
\begin{subequations}\label{eq:linODEs}
\begin{align}
\dot p_r&=\DpCoeff{}(q_r-q_l)\label{eq:dprdt_iso}\\
\dot q_\ell&=
\DqCoeffpr{} p_r+\DqCoeffpl{} p_\ell+\DqCoeffql{} q_\ell,\label{eq:dqldt_iso}
\end{align}
\end{subequations}
or equivalently,
\begin{subequations}
\begin{align*}
\dot x_{t} &= \begin{bmatrix}
0 & -\DpCoeff{}\\
\DqCoeffpr{} & \DqCoeffql{}
\end{bmatrix}x_{t}+\begin{bmatrix}
0 & \DpCoeff{}\\
\DqCoeffpl{} & 0
\end{bmatrix}u_{t},\\
y_{t} &= x_{t},
\end{align*}
\end{subequations}
with $x_{t}=\begin{bmatrix}
\tilde p_{r}&\tilde q_{\ell}
\end{bmatrix}^\top$ as the state vector and $u_{t}=\begin{bmatrix}
 \tilde p_{\ell} &\tilde q_{r}
\end{bmatrix}^\top$ as the input vector.
The coefficients are
\begin{align*}
\DpCoeff{}&=-\frac{R_sT_0z_0}{AL},\quad \DqCoeffpr{}=-\frac{A}{L},\\
 \DqCoeffpl{}&=\frac{A}{L}+\frac{\lambda R_sT_0z_0}{2DA}\frac{ q_{ss}\vert  q_{ss}\vert}{ p_{\ell,ss}^2}-\frac{Agh}{R_sT_0z_0L},\\
\DqCoeffql{}&=-\frac{\lambda R_sT_0z_0}{DA} \frac{|q_{ss}|}{ p_{\ell, ss}}, 
\end{align*}
for which the parameters are described in Table \ref{tab:parameters}. For clarity and without loss of generality, we assume throughout this work that the nominal mass flow $ q_{ss}$ is positive, although we stress that the denomination of the presented models corresponds to their positive $x$ direction and not necessarily to their flow direction. For instance, the joint introduced below could be a geometry where the flow physically branches.
\begin{table}[ht]
\begin{center}
\begin{tabular}{|c|l|l|}
\hline
Symbol&Meaning&SI-unit\\
\hline\hline
$A$&Cross-sectional area&$\scriptstyle[\text{m}^2]$\\
\hline
$D$&Pipe inside diameter&$\scriptstyle[\text{m}]$\\
\hline
$g$&Gravity constant&$\scriptstyle[\frac{\text{m}}{\text{s}^2}]$\\
\hline
$h(x)$&Pipe elevation from $x=0$ to $x=L$&$\scriptstyle[\text{m}]$\\
\hline
$L$&Pipe length&$\scriptstyle[\text{m}]$\\
\hline
$p(x,t)$&Pressure&$\scriptstyle[\frac{\text{kg}}{\text{s}^2\text{m}}]$\\
\hline
$q(x,t)$&Mass flow&$\scriptstyle[\frac{\text{kg}}{\text{s}}]$\\
\hline
$R_s$&Specific gas constant&$\scriptstyle[\frac{\text{m}^2}{\text{s}^2\text{K}}]$\\
\hline
$T_0$&Constant temperature&$\scriptstyle[\text{K}]$\\
\hline
$z_0$&Constant compressibility factor&$\scriptstyle[1]$\\
\hline
$\lambda$&Friction factor&$\scriptstyle[1]$\\
\hline
$(\cdot)_{ss}$&Nominal value&\\
\hline
\end{tabular}
\end{center}
\caption{Definitions of variables and SI-units.}
\label{tab:parameters}
\end{table}
\section{Pipes in series}
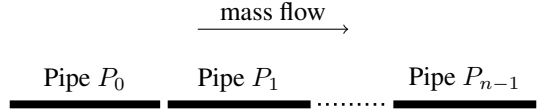
\begin{figure}[ht!]
    \centering
   \begin{tikzpicture}

\draw [->] (2.5,0) -- node[above]{mass flow} (4.5,0);
\draw [line width=3pt] (0,-1) -- node[above] {Pipe $P_{0}$}  (2,-1);
\draw [line width=3pt] (2.1,-1) -- node[above] {Pipe $P_{1}$}  (4,-1);
\draw[dotted, line width=1pt] (4.1,-1) -- (5,-1);
\draw [line width=3pt] (5.1,-1) -- node[above] {Pipe $P_{n-1}$}  (7,-1);
%
%
\end{tikzpicture}
 \caption{Series of $n$ pipes.}
    \label{fig:series_sketch}
\end{figure}

The single-pipe model from above is used to generate a composite model for multiple pipes in series as shown in Figure \ref{fig:series_sketch}, subsuming intermediate boundary conditions making them well-suited candidates for model-based control design.
The state, input and output elements are composed in lexicographical order, i.e., 
\begin{align*}
x_t&=\begin{bmatrix}
p_{0,r} & \dots & p_{n-1,r} & q_{0,\ell}&\dots & q_{n-1,\ell}
\end{bmatrix}^\top,\\
 u_t&=\begin{bmatrix}
p_{0,\ell} & q_{n-1,r}
\end{bmatrix}^\top, \\
y_t&=\begin{bmatrix}
p_{n-1,r}&q_{0,\ell}
\end{bmatrix}^\top.
\end{align*}
They yield \cite{sven_bob_rob_gas} 
\begin{align}\label{eq:multi_pipe_sys}
\dot x_t=A_sx_t+B_su_t,\qquad \qquad y_t=C_sx_t,
\end{align}
where
\begin{align*}
A_s&=\begin{bmatrix}
\zero & A_{s,12}\\
A_{s,21} & A_{s,22}
\end{bmatrix},\quad B_s=\begin{bmatrix}
B_{s,1}^\top & B_{s,2}^\top
\end{bmatrix}^\top,\quad
C_s=\begin{bmatrix}
\zero_{2,2(n-1)}&\eye_2&\zero_{2,2(n-1)}
\end{bmatrix},
\end{align*}
with $\zero$ being the zero matrix of appropriate size or, if appropriate, size indicated by the subscript. Accordingly, $\eye$ represents the identity matrix. Further, 
\begin{align*}
A_{s,12}&=\begin{bmatrix}
-\DpCoeff{0} & \DpCoeff{0} &&\zero\\
&\ddots & \ddots&\\
&& \ddots&\DpCoeff{n-2}\\
\zero&& &-\DpCoeff{n-1}\\
\end{bmatrix},\quad
A_{s,21}&=\begin{bmatrix}
\DqCoeffpr{0} &&&\zero\\
\DqCoeffpl{1}&\ddots & &\\
&\ddots& \ddots&\\
\zero&& \DqCoeffpl{n-1}&\DqCoeffpr{n-1}\\
\end{bmatrix}, \quad
A_{s,22}=\diag(\DqCoeffql{0}, \DqCoeffql{1}, \dots, \DqCoeffql{n-1}),\\
B_{s,1}&=\begin{bmatrix}
\zero_{2(n-1),2} \\ \begin{matrix}
0 & \DpCoeff{n-1}
\end{matrix}
\end{bmatrix},\quad B_{s,2}=\begin{bmatrix}
\begin{matrix}
\DqCoeffpl{0} & 0
\end{matrix}\\ \zero_{2(n-1),2}
\end{bmatrix}.
\end{align*}
\subsection*{Example for up to three pipes in series: frequency response}
\begin{figure}[ht!]
    \centering
    \includegraphics[scale=.55]{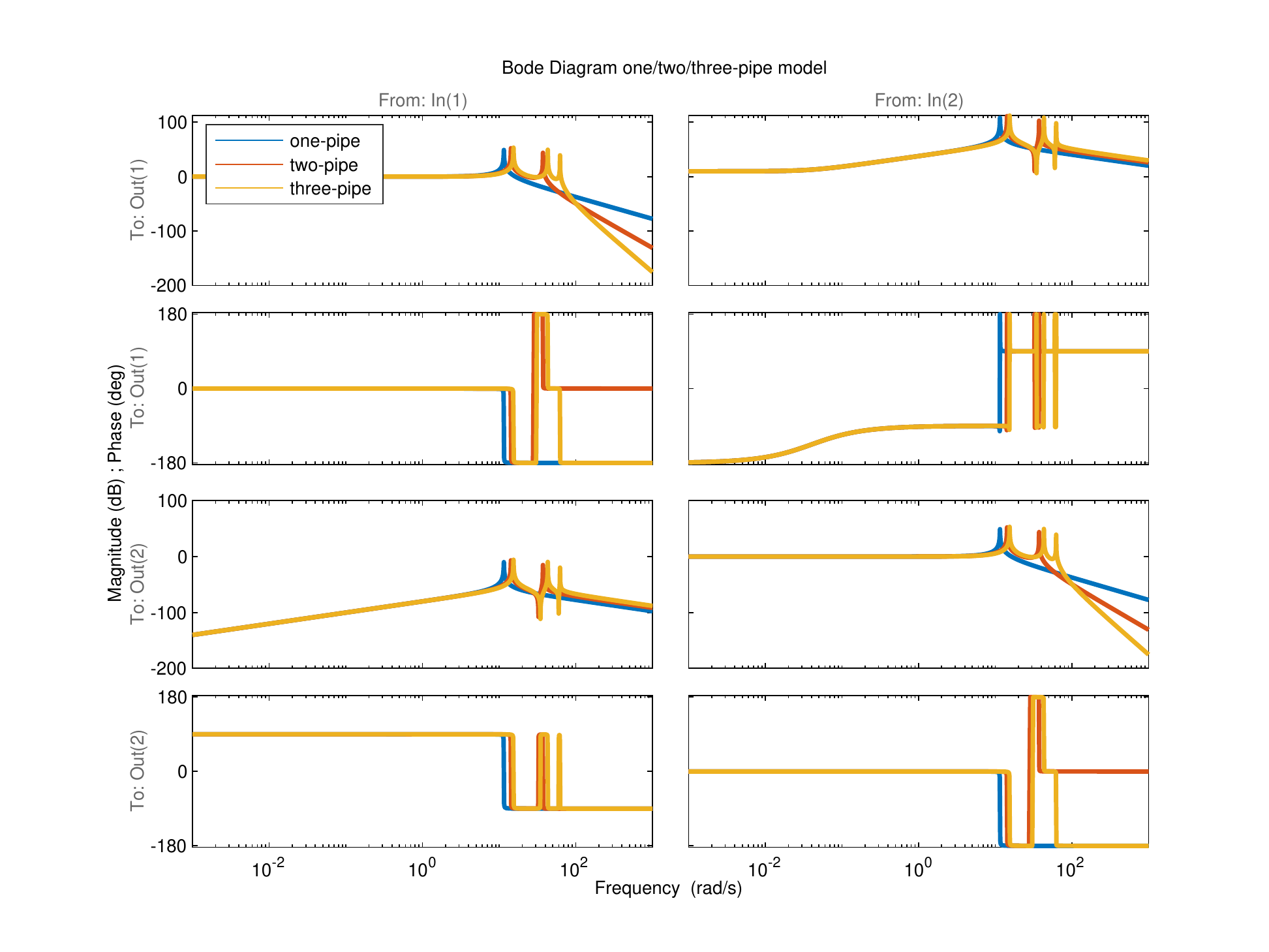}
 \caption{Comparison of the frequency responses from $\{p_{0,\ell},q_{n-1,r}\}\to\{p_{n-1,r},q_{0,\ell}\}$ with $n\in\{1,2,3\}$ between one pipe ($n=1$), two pipes $(n=2)$ and three pipes $(n=3$) in series with overall identical length.}
    \label{fig:123pipeModel}
\end{figure}

Next, in Figure \ref{fig:123pipeModel} we compare the frequency response between a pipe section of length $x$m modeled as: a single $x$m-long pipe using \eqref{eq:multi_pipe_sys} with $n=1$, two $x/2$m-long pipes in series using \ref{eq:multi_pipe_sys} with $n=2$, and three $x/3$m-long pipes in series using \eqref{eq:multi_pipe_sys} with $n=3$. 
We observe that the frequency responses at low frequencies coincide: conservation of mass is captured due to $0$dB gain from mass flow to mass flow in row three and column two; the mass flow shows derivative behavior with respect to pressure for row three/four and column one, an increase in mass flow leads to a drop in pressure in row one/two and column two as the phase is $-180^\circ$ for $\omega\to 0$. Differences only occur for high frequencies that lie outside our area of interest for slow process control.
\section{Joint}
We wish to model a pipe joint where several pipes merge into one, see sketch below in Figure \ref{fig:xjoiningPipes_sketch}. This extends \cite{sven_bob_rob_gas} where a joint of only two pipes is presented.
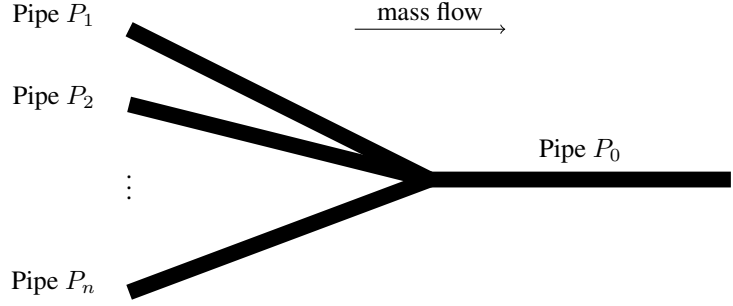
\begin{figure}[t!]
    \centering
   \begin{tikzpicture}

\draw [->] (3,0) -- node[above]{mass flow} (5,0);
\draw [line width=6pt] (0,0) -- node[above, xshift=-3cm, yshift=0.8cm] {Pipe $P_1$}  (4,-2);
\draw [line width=6pt] (0,-1) -- node[above, xshift=-3cm, yshift=0.2cm] {Pipe $P_2$}  (4,-2);
\node [line width=6pt] at (0,-2) {$\vdots$};
\draw [line width=6pt] (0,-3.5) -- node[above, xshift=-3cm, yshift=-1cm] {Pipe $P_n$}  (4,-2);
\draw [line width=6pt] (4,-2) -- node[above] {Pipe $P_{0}$}  (8,-2);
%
%
\end{tikzpicture}
 \caption{Joint of $n$ pipes merging into one}
    \label{fig:xjoiningPipes_sketch}
\end{figure}
Building on \eqref{eq:linODEs}, the joining pipes can be described by
\begin{align}
\dot p_{j, r}&=\DpCoeff{j} (q_{j,r}-q_{j,\ell})\label{eq:pdot},\\
\dot q_{j,\ell}&=\DqCoeffpr{j} p_{j,r}+\DqCoeffpl{j} p_{j,\ell}+\DqCoeffql{j} q_{j,\ell}\label{eq:qdot},
\end{align}
where the subscript $j$ corresponds to pipe $P_j$, with coefficients as defined above and $j\in\{0,1,\dots,n\}$. 
We say that each pair $(p_{j,r},q_{j,\ell})$ is the state variable and $p_{j,\ell}$ and $q_{j,r}$ are the input variables related to pipe $P_j$. This is consistent with the state-space realization above. 

\begin{defin}
 We say \textup{total inputs} for inputs into the joint model that are not related to states of any other pipes via algebraic constraints. We say \textup{internal variables} for input variables that are not total inputs. 
\end{defin} 
 Subsequently, towards induction, we derive state-space models for $n=2,3$ first and then extend the result to $n$ elements. 
\subsection{Two joining pipes}
Consider algebraic constraints
\begin{align}
p_{1,r}=p_{2,r}=p_{0,\ell},\label{eq:ACpn=2}\\
q_{1,r}+q_{2,r}=q_{0,\ell},\label{eq:ACqn=2}
\end{align}
which express conservation of mass and continuity of the pressure.
To eliminate internal variables $q_{1,r},q_{r,2}$, consider \eqref{eq:ACpn=2}, take the derivative and use \eqref{eq:pdot} to obtain
\begin{align}
(\DpCoeff{1}+\DpCoeff{2})q_{1,r}&=\DpCoeff{1}q_{1,\ell}+\DpCoeff{2}(q_{0,\ell}-q_{2,\ell}),\label{eq:q1r_n2}\\
(\DpCoeff{1}+\DpCoeff{2})q_{2,r}&=\DpCoeff{2}q_{2,\ell}+\DpCoeff{1}(q_{0,\ell}-q_{1, \ell})\label{eq:q2r_n2},
\end{align}
which gives rise to a state formulation of
\begin{align*}
\dot{\begin{bmatrix}
p_{0,r}\\q_{0,\ell}\\p_{1,r}\\q_{1,\ell}\\q_{2,\ell}
\end{bmatrix}}&=
\begin{bmatrix}
0&-\DpCoeff{0}&0&0&0\\
\DqCoeffpr{0}&\DqCoeffql{0}& \DqCoeffpl{0}&0&0\\
0&\frac{\DpCoeff{1}\DpCoeff{2}}{\DpCoeff{1}+\DpCoeff{2}}&0 & \frac{\DpCoeff{1}^2}{\DpCoeff{1}+\DpCoeff{2}}-\DpCoeff{1}&-\frac{\DpCoeff{1}\DpCoeff{2}}{\DpCoeff{1}+\DpCoeff{2}}\\
0&0&\DqCoeffpr{1}&\DqCoeffql{1}&0\\
0 & 0 & \DqCoeffpr{2} & 0 & \DqCoeffql{2}
\end{bmatrix}
\begin{bmatrix}
p_{0,r}\\q_{0,\ell}\\p_{1,r}\\q_{1,\ell}\\q_{2,\ell}
\end{bmatrix}+
\begin{bmatrix}
0 & 0 & \DpCoeff{0}\\
0 & 0 & 0\\
0 & 0 & 0\\
\DqCoeffpl{1} & 0 & 0\\
0 & \DqCoeffpl{2} & 0
\end{bmatrix}
\begin{bmatrix}
p_{1,\ell}\\p_{2,\ell}\\q_{0,r}
\end{bmatrix},\\
\begin{bmatrix}
p_{0,r}\\q_{1,\ell}\\q_{2,\ell}
\end{bmatrix}&=\begin{bmatrix}
1 & 0 & 0 & 0 & 0\\
0 & 0 & 0 & 1 & 0\\
0 & 0 & 0 & 0 & 1\\
\end{bmatrix}
\begin{bmatrix}
p_{0,r}\\q_{0,\ell}\\p_{1,r}\\q_{1,\ell}\\q_{2,\ell}
\end{bmatrix}.
\end{align*}
We wish to extend this result to any finite $n$. In particular, we would like to obtain an expression for $q_{j,r}$ only depending on state variables and total inputs, as in \eqref{eq:q1r_n2} and \eqref{eq:q2r_n2}. Toward this goal, we derive expressions for $q_{r, i},\,i\in\{1,2,3\},$ for $n=3$.

\subsection{Three joining pipes}
Consider the algebraic constraints for $n=3$,
\begin{align*}
p_{1,r}=p_{2,r}=p_{3,r}=p_{0,\ell},\\
q_{1,r}+q_{2,r}+q_{3,r}=q_{0,\ell}.
\end{align*}
Then,
\begin{align*}
\dot p_{1,r}&=\dot p_{2,r},\\
(\DpCoeff{1}+\DpCoeff{2})q_{1,r}&=\DpCoeff{1}q_{1,\ell}+\DpCoeff{2}(q_{0,\ell}-q_{2,\ell}-q_{3,r}),\\[1em]
\dot p_{3,r}&=\dot p_{2,r},\\
(\DpCoeff{2}+\DpCoeff{3})q_{3,r}&=\DpCoeff{3}q_{3,\ell}+\DpCoeff{2}(q_{0,\ell}-q_{2,\ell}-q_{1,r})\\
&=\DpCoeff{3}q_{3,\ell}+\DpCoeff{2}\Big(q_{0,\ell}-q_{2,\ell}-(\DpCoeff{1}+\DpCoeff{2})^{-1}[\DpCoeff{1}q_{1,\ell}+\DpCoeff{2}(q_{0,\ell}-q_{2,\ell}-q_{3,r})]\Big).
\end{align*}
Note that 
\begin{align*}
q_{3,r}: \quad &(\DpCoeff{2}+\DpCoeff{3})-\frac{\DpCoeff{2}^2}{\DpCoeff{1}+\DpCoeff{2}}=\frac{\DpCoeff{1}\DpCoeff{2}+\DpCoeff{1}\DpCoeff{3}+\DpCoeff{2}\DpCoeff{3}}{\DpCoeff{1}+\DpCoeff{2}},\\
q_{0, \ell} ,\,(q_{2,\ell}): \quad & \DpCoeff{2}-\frac{\DpCoeff{2}^2}{\DpCoeff{1}+\DpCoeff{2}}=\frac{\DpCoeff{1}\DpCoeff{2}}{\DpCoeff{1}+\DpCoeff{2}} \left(-\frac{\DpCoeff{1}\DpCoeff{2}}{\DpCoeff{1}+\DpCoeff{2}}\right),\\
q_{3,\ell}:\quad& \DpCoeff{3} = \frac{\DpCoeff{1}\DpCoeff{3}+\DpCoeff{2}\DpCoeff{3}}{\DpCoeff{1}+\DpCoeff{2}},\\
q_{1,\ell}:\quad& -\frac{\DpCoeff{1}\DpCoeff{2}}{\DpCoeff{1}+\DpCoeff{2}},
\end{align*}
so that
\begin{align}
(\DpCoeff{1}\DpCoeff{2}+\DpCoeff{1}\DpCoeff{3}+\DpCoeff{2}\DpCoeff{3})q_{3,r}=\DpCoeff{1}\DpCoeff{2}(q_{0,\ell}-q_{1,\ell}-q_{2,\ell})+(\DpCoeff{1}\DpCoeff{3}+\DpCoeff{2}\DpCoeff{3})q_{3,\ell}.\label{eq:q3r_n3}
\end{align}
Similarly,
\begin{align}
(\DpCoeff{1}\DpCoeff{2}+\DpCoeff{1}\DpCoeff{3}+\DpCoeff{2}\DpCoeff{3})q_{1,r}=\DpCoeff{2}\DpCoeff{3}(q_{0,\ell}-q_{2,\ell}-q_{3,\ell})+(\DpCoeff{1}\DpCoeff{2}+\DpCoeff{1}\DpCoeff{3})q_{1,\ell},\label{eq:q1r_n3}\\
(\DpCoeff{1}\DpCoeff{2}+\DpCoeff{1}\DpCoeff{3}+\DpCoeff{2}\DpCoeff{3})q_{2,r}=\DpCoeff{1}\DpCoeff{3}(q_{0,\ell}-q_{1,\ell}-q_{3,\ell})+(\DpCoeff{1}\DpCoeff{2}+\DpCoeff{2}\DpCoeff{3})q_{2,\ell}.\label{eq:q2r_n3}
\end{align}
Comparing equations \eqref{eq:q3r_n3}--\eqref{eq:q2r_n3} with \eqref{eq:q1r_n2} and \eqref{eq:q2r_n2} we recognize a certain pattern which we explore next.
\subsection{$\mathbf{n}$ joining pipes}
Facilitated by observations above we derive a general statement.
\begin{theorem}\label{thm:qkr}
Consider $n$ joining pipes with individual dynamics \eqref{eq:pdot} - \eqref{eq:qdot}. Let $q_{j,\ell},\, j\in\{0,1,\dots,n\},$ be state variables. The internal variables $q_{k,r},\,k\in\{1,2,\dots,n\},$ can be described by state variables,
\begin{align}
\left(\sum_{j=1}^n\prod_{\substack{i=1\\i\neq j}}^n\DpCoeff{i}\right)q_{k,r}=\prod_{\substack{i=1\\i\neq k}}^n \DpCoeff{i}\left(q_{0,\ell}-\sum_{\substack{i=1\\i\neq k}}^n q_{i,\ell}\right)+\left(\sum_{j=1}^n\prod_{\substack{i=1\\i\neq j}}^n\DpCoeff{i}-\prod_{\substack{i=1\\i\neq k}}^n\DpCoeff{i}\right)q_{k,\ell}\label{eq:qkr}
\end{align}
\end{theorem}
\begin{proof}
We prove the theorem by induction. We observe that \eqref{eq:qkr} holds for $n\in\{2,3\}$. Hence, assume it holds for some $n$ and for $n+1$ note that interconnections dictate the algebraic constraints 
\begin{align*}
p_{1,r}=p_{2,r}=\dots=&p_{n+1,r}=p_{0,\ell},\\
q_{1,r}+q_{2,3}+\dots+q_{r,n}+q_{r,n+1}&\doteq q^n_{0,\ell}+q_{r,n+1}=q_{0,\ell}.
\end{align*}
Hence, by the algebraic constraints on the pressure and the pipe dynamics, 
\begin{align*}
q_{n+1,r}&=\frac{1}{\DpCoeff{n+1}}(\dot p_{k,r}-\dot p_{n+1,r})+q_{n+1,r}\\
&=\frac{1}{\DpCoeff{n+1}}\left(\DpCoeff{k}(q_{k,r}-q_{k,\ell})-\dot p_{n+1,r}\right)+q_{n+1,\ell}+\frac{\dot p_{n+1,r}}{\DpCoeff{n+1}}\\
&=\frac{1}{\DpCoeff{n+1}}\DpCoeff{k}(q_{k,r}-q_{k,\ell})+q_{n+1,\ell},
\end{align*}
where without loss of generality $k\in\{1,2\dots,n-1\}$.
Then \eqref{eq:qkr} yields
\begin{align*}
\left(\sum_{j=1}^n\prod_{\substack{i=1\\i\neq j}}^n\DpCoeff{i}\right)q_{k,r}&=\prod_{\substack{i=1\\i\neq k}}^n \DpCoeff{i}\left(q^n_{0,\ell}-\sum_{\substack{i=1\\i\neq k}}^n q_{i,\ell}\right)+\left(\sum_{j=1}^n\prod_{\substack{i=1\\i\neq j}}^n\DpCoeff{i}-\prod_{\substack{i=1\\i\neq k}}^n\DpCoeff{i}\right)q_{k,\ell}\\
&=\prod_{\substack{i=1\\i\neq k}}^n \DpCoeff{i}\left(q_{0,\ell}-q_{n+1,r}-\sum_{\substack{i=1\\i\neq k}}^n q_{i,\ell}\right)+\left(\sum_{j=1}^n\prod_{\substack{i=1\\i\neq j}}^n\DpCoeff{i}-\prod_{\substack{i=1\\i\neq k}}^n\DpCoeff{i}\right)q_{k,\ell}\\
&=-\prod_{\substack{i=1}}^n \DpCoeff{i}\frac{1}{\DpCoeff{n+1}}(q_{k,r}-q_{k,\ell})+ \prod_{\substack{i=1\\i\neq k}}^n \DpCoeff{i}\left(q_{0,\ell}-q_{n+1,\ell}-\sum_{\substack{i=1\\i\neq k}}^n q_{i,\ell}\right)+\left(\sum_{j=1}^n\prod_{\substack{i=1\\i\neq j}}^n\DpCoeff{i}-\prod_{\substack{i=1\\i\neq k}}^n\DpCoeff{i}\right)q_{k,\ell}\\
&=-\frac{1}{\DpCoeff{n+1}}\prod_{\substack{i=1}}^n \DpCoeff{i}q_{k,r}+ \prod_{\substack{i=1\\i\neq k}}^n \DpCoeff{i}\left(q_{0,\ell}-\sum_{\substack{i=1\\i\neq k}}^{n+1} q_{i,\ell}\right)+\left(\sum_{j=1}^n\prod_{\substack{i=1\\i\neq j}}^n\DpCoeff{i}-\prod_{\substack{i=1\\i\neq k}}^n\DpCoeff{i}+\frac{1}{\DpCoeff{n+1}}\prod_{\substack{i=1}}^n \DpCoeff{i}\right)q_{k,\ell}.
\end{align*}
Move the first term of the right-hand side to the left, multiply the equation by $\DpCoeff{n+1}$ and observe that the left-hand side is
\begin{align*}
\sum_{j=1}^n\prod_{\substack{i=1\\i\neq j}}^n\DpCoeff{i}\DpCoeff{n+1}+\prod_{\substack{i=1}}^n \DpCoeff{i}=\sum_{j=1}^n\prod_{\substack{i=1\\i\neq j}}^{n+1}\DpCoeff{i}+\prod_{\substack{i=1}}^n \DpCoeff{i}=\sum_{j=1}^{n+1}\prod_{\substack{i=1\\i\neq j}}^{n+1}\DpCoeff{i}.
\end{align*}
This identity can also be used for the last term on the right-hand side. Then,
\begin{align*}
\left(\sum_{j=1}^{n+1}\prod_{\substack{i=1\\i\neq j}}^{n+1}\DpCoeff{i}\right)q_{k,r}&=\prod_{\substack{i=1\\i\neq k}}^{n+1} \DpCoeff{i}\left(q_{0,\ell}-\sum_{\substack{i=1\\i\neq k}}^{n+1} q_{i,\ell}\right)+ \left(\sum_{j=1}^{n+1}\prod_{\substack{i=1\\i\neq j}}^{n+1}\DpCoeff{i}-\prod_{\substack{i=1\\i\neq k}}^{n+1}\DpCoeff{i}\right)q_{k,\ell},
\end{align*}
which is equivalent to our induction hypothesis and thus concludes the proof.
\end{proof}
Theorem \ref{thm:qkr} enables us to formulate a general state-space model for a joint of $n$ pipes for which we define $\zero_i$ $(\zero_{i,j})$ as a vector (matrix) of zeros of dimension $i$ ($i\times j$). The vector and matrix of ones $\one_i$ and $\one_{i,j}$ are denoted accordingly.
\begin{corollary}\label{cor:state-space-njoint}
A joint of $n$ pipes with dynamics \eqref{eq:pdot} - \eqref{eq:qdot} can be described by the state-space model,
\begin{align*}
\dot{x}&=Ax+Bu=\begin{bmatrix}
\zero_{2,2} & A_{12}\\
A_{21} & A_{22}
\end{bmatrix} x+\begin{bmatrix}
B_1\\
B_2
\end{bmatrix} u,\\
y&=Cx=\begin{bmatrix}
\begin{matrix}
1 & 0 & 0 & 0\\
0 & 0 & 0 & 1
\end{matrix} & \zero_{2,n-1}\\\zero_{n-1,4}&\eye_{n-1}
\end{bmatrix} x,
\end{align*}
where
\begin{align*}
u&=
\begin{bmatrix}
p_{1,\ell}&p_{2,\ell}&\dots&p_{n,\ell}&q_{0,r}
\end{bmatrix}^\top\in\R^{n+1},\\
x&=
\begin{bmatrix}
p_{0,r}&p_{1,r}&q_{0,\ell}&q_{1,\ell}&q_{2,\ell}&\dots&q_{n,\ell}
\end{bmatrix}^\top\in\R^{n+3}, \\
y&=
\begin{bmatrix}
p_{0,r}&q_{1,\ell}&q_{2,\ell}&\dots&q_{n,\ell}
\end{bmatrix}^\top\in\R^{n+1},
\end{align*}
\begin{align*}
A_{12}&=
\begin{bmatrix}
-\DpCoeff{0} & \zero_{n}^\top\\
a & -a\one_{n}^\top
\end{bmatrix},\quad
A_{21}=
\begin{bmatrix}
\DqCoeffpr{0} & 0 & 0 & \dots & 0\\
\DqCoeffpl{0} &  \DqCoeffpr{1} &  \DqCoeffpr{2}& \dots & \DqCoeffpr{n}
\end{bmatrix}^\top,\quad
A_{22}=\diag(\DqCoeffql{0} ,\dots,\DqCoeffql{n}) \\
B_1&=
\begin{bmatrix}
\zero_n^\top& \DpCoeff{0} \\
\zero_n^\top& 0
\end{bmatrix},\quad 
B_2=\begin{bmatrix}
\zero^\top_{n+1}\\
\begin{matrix}
\begin{matrix}
\DqCoeffpl{1}&&\\
&\ddots&\\
&&\DqCoeffpl{n}
\end{matrix}&\zero_n
\end{matrix}
\end{bmatrix}
\end{align*}
with
\begin{align*}
a&=\DpCoeff{1}\left(\sum_{j=1}^n\prod_{\substack{i=1\\i\neq j}}^n\DpCoeff{i}\right)^{-1}\prod_{\substack{i=2}}^n \DpCoeff{i}.
\end{align*}
\end{corollary}
\begin{proof}
The state space realization is a direct consequence of the pipe dynamics in \eqref{eq:pdot}-\eqref{eq:qdot}, the algebraic constraint of equal pressure at the intersection, Theorem \ref{thm:qkr} and the observation that for the ODE of $p_{1,r}$,
\begin{align*}
%
\DpCoeff{1}\left(\sum_{j=1}^n\prod_{\substack{i=1\\i\neq j}}^n\DpCoeff{i}\right)^{-1}\left(\sum_{j=1}^n\prod_{\substack{i=1\\i\neq j}}^n\DpCoeff{i}-\prod_{\substack{i=2}}^n\DpCoeff{i}\right)-\DpCoeff{1}
%
&=-\DpCoeff{1}\left(\sum_{j=1}^n\prod_{\substack{i=1\\i\neq j}}^n\DpCoeff{i}\right)^{-1}\prod_{\substack{i=2}}^n \DpCoeff{i}.
\end{align*} 
\end{proof}
As a consequence of the algebraic constraints, the state-space realization above directly captures conservation of mass.
\begin{corollary}\label{cor:conservationMass-state-space-njoint}
The state-space realization in Corollary \ref{cor:state-space-njoint} satisfies conservation of mass, i.e., at steady state,
\begin{align*}
q_{0,\ell}=q_{1,r}+q_{2,r}+\dots+q_{n,r}.
\end{align*}
\end{corollary}
\begin{proof}
The second row of the steady state equation $\zero_{n+3}=Ax+Bu$ yields the desired result.
\end{proof}

\section{Star junction}
\begin{figure}[ht!]
    \centering
   \begin{tikzpicture}

\draw [->] (3,0) -- node[above]{mass flow} (5,0);
\draw [line width=6pt] (0,0) -- node[above, xshift=-3cm, yshift=0.8cm] {Pipe $P_1$}  (4,-2);
\draw [line width=6pt] (0,-1) -- node[above, xshift=-3cm, yshift=0.2cm] {Pipe $P_2$}  (4,-2);
\node [line width=6pt] at (0,-2) {$\vdots$};
\draw [line width=6pt] (0,-3.5) -- node[above, xshift=-3cm, yshift=-1cm] {Pipe $P_n$}  (4,-2);
\draw [line width=6pt] (4,-2) -- node[above, xshift=3cm, yshift=0.8cm] {Pipe $P_{n+1}$}  (8,0);
\draw [line width=6pt] (4,-2) -- node[above, xshift=3cm, yshift=0.2cm] {Pipe $P_{n+2}$}  (8,-1);
\node [line width=6pt] at (8,-2) {$\vdots$};
\draw [line width=6pt] (4,-2) -- node[above, xshift=3cm, yshift=-1cm] {Pipe $P_{n+m}$}  (8,-3.5);
%
%
\end{tikzpicture}
 \caption{Star junction of $n+m$ pipes}
    \label{fig:star_junction_sketch}
\end{figure}
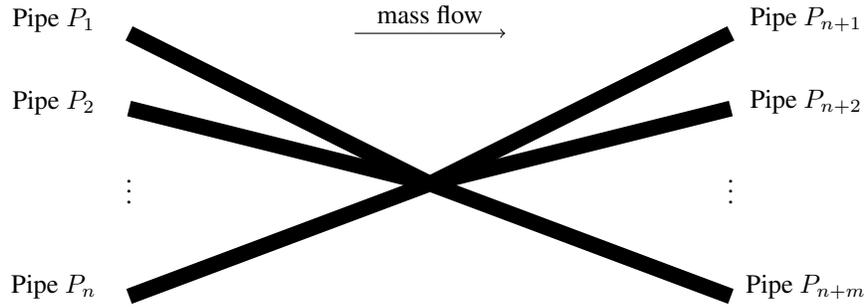
Theorem \ref{thm:qkr} also facilitates a statement about a general star junction, which consists of $n$ pipes connecting to $m$ pipes, as illustrated in Figure \ref{fig:star_junction_sketch}. Towards a state-space realization and akin to the case of a joint, we wish to express internal variables in terms of state variables.
\begin{corollary}\label{col:pkr_qkr}
Consider a star junction with $n$ joining pipes connecting to $m$ branching pipes with individual dynamics \eqref{eq:pdot} - \eqref{eq:qdot}. Let $q_{i,\ell},\, i\in\{1,2,\dots,n+m\}$, and $p_{1,r}$ be state variables. The internal variables $q_{k,r}, \,k\in\{1,2,\dots,n\}$ and $p_{j,\ell},\,j\in\{n+1,n+2,\dots,n+m\},$ can be described by state variables through 
\begin{align}
\left(\sum_{j=1}^n\prod_{\substack{i=1\\i\neq j}}^n\DpCoeff{i}\right)q_{k,r}&=\prod_{\substack{i=1\\i\neq k}}^n \DpCoeff{i}\left(\sum_{i=n+1}^{n+m}q_{i,\ell}-\sum_{\substack{i=1\\i\neq k}}^n q_{i,\ell}\right)+\left(\sum_{j=1}^n\prod_{\substack{i=1\\i\neq j}}^n\DpCoeff{i}-\prod_{\substack{i=1\\i\neq k}}^n\DpCoeff{i}\right)q_{k,\ell}\label{eq:qkr_star_junction},\\
p_{j,\ell}&=p_{1,r}\label{eq:pjl_star}.
\end{align}
\end{corollary}
\begin{proof}
Consider \eqref{eq:qkr} and note that therein that $q_{0,\ell}=\sum_{j=n+1}^{n+m}q_{j,\ell}$. The result for the first equation follows immediately. The second equation is a direct consequence of the algebraic constraint of all pressures being equal at the intersection.
\end{proof}
We are now able to derive a state-space realization.
\begin{corollary}
A star junction of $n$ pipes connecting to $m$ pipes with individual dynamics \eqref{eq:pdot} - \eqref{eq:qdot} can be represented by the state-space model,
\begin{align*}
\dot x_1=A_1x_1+B_1u, \quad \dot x_2&=A_1x_2+B_2u, \quad \dot x_3=A_1x_3+B_3u,\quad \dot x_4=A_1x_4+B_4u,\\
y&=C(x_1+x_2+x_3+x_4),
\end{align*}
where
\begin{align*}
x_1&=x_2=x_3=x_4=
\begin{bmatrix}
p_{1,r}\\p_{n+1,r}\\p_{n+2,r}\\\vdots\\p_{n+m,r}\\q_{1,\ell}\\q_{2,\ell}\\\vdots\\q_{n+m,\ell}
\end{bmatrix}\in\R^{n+2m+1}, \quad 
u=
\begin{bmatrix}
p_{1,\ell}\\p_{2,\ell}\\\vdots\\p_{n,\ell}\\q_{n+1,r}\\q_{n+2,r}\\\vdots\\q_{n+m,r}
\end{bmatrix}\in\R^{n+m},\quad
y=
\begin{bmatrix}
p_{n+1,r}\\p_{n+2,r}\\\vdots\\p_{n+m,r}\\q_{1,\ell}\\q_{2,\ell}\\\vdots\\q_{n,\ell}
\end{bmatrix}\in\R^{n+m},\\
A_1&=
\begin{bmatrix}
\begin{array}{c}
\begin{matrix}
\mathbf{0}_{m+1}^\top & -a & -a\mathbf{1}_{n-1}^\top & a\mathbf{1}_{m}^\top
\end{matrix}\\
\hline
\mathbf{0}_{n+2m,n+2m+1}
\end{array}
\end{bmatrix},
\quad
B_1=\mathbf{0}_{n+2m+1,n+m},\\
A_2 &= \begin{bmatrix}
\begin{array}{c|c}
\mathbf{0}_{n+2m+1,n+m+1} &
\begin{array}{c}
\begin{matrix}
\mathbf{0}_{m}^\top\\
\hline
\begin{matrix}
-\DpCoeff{n+1}&&\mathbf{0}\\
&\ddots&\\
\mathbf{0}&&-\DpCoeff{n+m}
\end{matrix}
\end{matrix}\\
\hline
\mathbf{0}_{n+m}^\top
\end{array} 
\end{array}
\end{bmatrix},\quad
B_2=
\begin{bmatrix}
\begin{array}{c}
\mathbf{0}_{n+m}^\top\\
\hline
\begin{array}{c|c}
\mathbf{0}_{m,n}
&
\begin{matrix}
\DpCoeff{n+1}&&\mathbf{0}\\
&\ddots&\\
\mathbf{0}&&\DpCoeff{n+m}
\end{matrix}
\end{array}\\
\hline
\mathbf{0}_{n+m,n+m}
\end{array}
\end{bmatrix},\\
A_3&=
\begin{bmatrix}
\begin{array}{c}
\mathbf{0}_{m+1,n+2m+1}\\
\hline
\begin{array}{c|c|c|c}
\begin{matrix}
\DqCoeffpr{1}\\
\vdots\\
\DqCoeffpr{n}
\end{matrix}
&
\mathbf{0}_{n,m}
&
\begin{matrix}
\DqCoeffql{1}&&\mathbf{0}\\
&\ddots&\\
\mathbf{0} &&\DqCoeffql{n}
\end{matrix}
&
\mathbf{0}_{n,m}
\end{array}\\
\hline
\mathbf{0}_{m,n+2m+1}
\end{array}
\end{bmatrix}, \quad
B_3=
\begin{bmatrix}
\begin{array}{c}
\mathbf{0}_{m+1,n+m}\\
\hline
\begin{array}{c|c}
\begin{matrix}
\DqCoeffpl{1}&&\mathbf{0}\\
&\ddots&\\
\mathbf{0}&&\DqCoeffpl{n}
\end{matrix}
&
\mathbf{0}_{n,m}
\end{array}\\
\hline
\mathbf{0}_{m,n+m}
\end{array}
\end{bmatrix},\\
A_4&=\begin{bmatrix}
\begin{array}{c}
\mathbf{0}_{n+m+1,n+2m+1}\\
\hline
\begin{array}{c|c|c|c}
\begin{matrix}
\DqCoeffpl{n+1}\\
\vdots\\
\DqCoeffpl{n+m}
\end{matrix}
&
\begin{matrix}
\DqCoeffpr{n+1}&&\mathbf{0}\\
&\ddots&\\
\mathbf{0} &&\DqCoeffpr{n+m}
\end{matrix}
&
\mathbf{0}_{m,n}
&
\begin{matrix}
\DqCoeffql{n+1}&&\mathbf{0}\\
&\ddots&\\
\mathbf{0} &&\DqCoeffql{n+m}
\end{matrix}
\end{array}
\end{array}
\end{bmatrix},\quad
B_4=\mathbf{0}_{n+2m+1,n+m},\\
C&=\begin{bmatrix}
\mathbf{0}_{n+m} & I_{n+m} & \mathbf{0}_{n+m,m}
\end{bmatrix},
\end{align*}
with
\begin{align*}
a&=\DpCoeff{1}\left(\sum_{j=1}^n\prod_{\substack{i=1\\i\neq j}}^n\DpCoeff{i}\right)^{-1}\prod_{\substack{i=2}}^n \DpCoeff{i}.
\end{align*}
\end{corollary}
\begin{proof}
The system related to $x_1$ describes the solution to $\dot p_{1,r}$ and directly follows from \eqref{eq:qkr_star_junction}. The system corresponding to $x_2$ describes $p_{j,r},\,j\in\{n+1,n+2,\dots,n+m\}$ and is immediately obtained by the dynamics \eqref{eq:pdot}. The systems for $x_3$ and $x_4$ relate to $q_{j,\ell},\,j\in\{1,2,\dots,n+m\}$ and are a consequence of the dynamics \eqref{eq:qdot} and the algebraic constraint \eqref{eq:pjl_star}.
\end{proof}

Corollary \ref{cor:conservationMass-state-space-njoint} showing conservation of mass at steady state can be extended to the star junction and is left as an exercise to the admittedly motivated reader.

\section{Conclusion}
We derived composite state-space models for the transient dynamics of one-dimensional isothermal gas flow through intersecting pipe geometries that are well-suited candidates for model-based control design. They also capture conservation of mass at steady state by subsuming algebraic constraints that would otherwise appear as part of a system of DAEs. Future research directions may include the analysis on when the given assumptions such as a constant temperature and the proposed algebraic constraints at the boundaries hold, especially related to the number of pipes, and a model validation using experimental data.

\section*{Acknowledgement}
This research was supported by funding from Solar Turbines {Incorporated}, who also provided operating data and guidance.

\bibliographystyle{plain}
\bibliography{bib_all.bib}

\end{document}